\documentstyle[11pt,titlepage,amssymb,euscript]{article}
\newtheorem{theorem}{Theorem}[section]
\newtheorem{proposition}[theorem]{Proposition}
\newtheorem{lemma}[theorem]{Lemma}
\newtheorem{corollary}[theorem]{Corollary}
\newtheorem{definition}[theorem]{Definition}
\newtheorem{example}[theorem]{Example}

\newcommand{\eop}{\hspace*{\fill}\mbox{ \lower.5ex\hbox{$\Box$}}}
\makeatletter
\def\proof{\@ifnextchar[{\@opargproof}{\@proof}}
\def\endproof{\eop\end{list}}

\def\@opargproof[#1]{\begin{list}{\bf Proof (#1):\hfill}%
{\parsep=\parskip\leftmargin=0pt\listparindent=\parindent
\settowidth{\labelwidth}{\bf Proof (#1):}\itemindent=\labelwidth
\addtolength{\itemindent}{\labelsep}}\item}
\let\end@opargproof=\endproof

\def\@proof{\begin{list}{\bf Proof:\hfill}%
{\parsep=\parskip\leftmargin=0pt\listparindent=\parindent
\settowidth{\labelwidth}{\bf Proof:}\itemindent=\labelwidth
\addtolength{\itemindent}{\labelsep}}\item}
\let\end@proof=\endproof

\@addtoreset{equation}{section}
\makeatother

\newcommand{\R}{\mbox{${\Bbb R}$}}

\newcommand{\g}{\mbox{$\gamma$}}

\newcommand{\ns}{\normalshape} 
\newcommand{\p}{\mbox{$\pi$}}
\newcommand{\scr}{\EuScript}

\renewcommand{\l}{\mbox{$\lambda$}}
\renewcommand{\sp}{\mbox{\ns Spray}}

\font\heads = cmbx12
\font\subheads = cmbx10 scaled 1095

\title{\Large\bf Pseudoconvex and Disprisoning\\Homogeneous Sprays}

\author{
 L. Del Riego\\
\noalign{\medskip}
{\small        Facultad de Ciencias}\\
{\small        Zona Universitaria}\\
{\small        Universidad Aut\'onoma de San Luis Potos\'\i}\\
{\small          San Luis Potos\'\i, SLP}\\
{\small               78290 MEXICO}\\
{\small     bestr1084@bestsd.sdsu.edu}\\
\and
Phillip. E. Parker\thanks{Partially supported
                                          by NSF grant INT-8921666.}\\
\noalign{\medskip}
{\small      Mathematics Department}\\
{\small            Wichita State University}\\
{\small             Wichita KS 67260-0033}\\
{\small                    USA}\\
{\small             pparker@twsuvm.uc.twsu.edu}
}

\date{
\vspace{.5in}
October 12, 1993\\
 \vspace{.5in}
 \begin{quote}
 {\bf Abstract}\ \
The pseudoconvex and disprisoning conditions for geodesics of linear
connections
are extended to the solution curves of general homogeneous sprays. The main
result is that pseudoconvexity and disprisonment are jointly stable in the fine
topology on the space of all homogeneous sprays of any degree of homogeneity.
 \end{quote}
\vspace{.4in}
MSC(1991) {\it primary\/} 53C15; {\it secondary\/} 53C22, 58F17.
}

\begin{document}

\maketitle

\section{\heads Introduction}
A spray is a natural generalization of the system of geodesics of a
linear connection. Pseudoconvexity and disprisonment are two increasingly
important properties which such a system may have. Since readers are not
likely to be equally familiar with all of these, we shall try in this
section to give enough background to provide some motivation.

\subsection{\subheads Pseudoconvexity and disprisonment}
Global hyperbolicity is well known to play an important role in
Lorentz\-ian
geometry and general relativity. It is a sufficient condition for the
existence of maximal length geodesic segments joining causally related
points \cite{S}. This gives a partial generalization of the important
Hopf-Rinow Theorem to Lorentzian manifolds. A spacetime has a global Cauchy
surface $N$ if and only if it is a globally hyperbolic spacetime \cite
{G}. Furthermore, these spacetimes are topological products of the form
$\R \times N$ and global hyperbolicity is a stable property. One
important
application  of global hyperbolicity is in the singularity theorems of
Hawking and Penrose \cite{H}. Many spacetimes have large globally
hyperbolic subsets, and these subsets may be used to construct causal
geodesics without conjugate points. The timelike convergence condition and
the generic condition then imply the incompleteness of such geodesics.
Existence of an incomplete causal geodesic is usually taken as indicating
a physical singularity.

Beem and Parker \cite{KGSGG,GBSS,WSS,PGR,PGC} have considered a
generalization of global hyperbolicity called {\em pseudoconvexity\/},
which
often can be used in place of global hyperbolicity. A spacetime $(M,g)$ is
said to be causally pseudoconvex if and only if given any compact set $K$
in $M$ there is always a larger compact set $K'$ such that all
causal
geodesics segments  joining points of $K$ lie entirely within $K'$.
This basic concept can be used for any class of geodesics. For example,
null pseudoconvexity is the requirement that all null geodesic segments
with endpoints in $K$ lie entirely within $K'$.

Like global hyperbolicity, (causal) pseudoconvexity is a type of
completeness requirement. Intuitively, one may think of pseudoconvex
spaces as failing to have any ``interior'' points missing. Thus, Minkowski
space less any compact set is neither globally hyperbolic nor causally
pseudoconvex.
A simple example of a causally pseudoconvex spacetime which is not
globally hyperbolic is the open strip $a < x < b$ in the Minkowski $(t,x)$
plane. We say that causal pseudoconvexity generalizes global
hyperbolicity since every globally hyperbolic spacetime is causally
pseudoconvex.

In the theory of pseudodifferential equations, the concept of
pseudoconvexity is applied to bicharacteristic segments in the study of
global solvability. If $(M,g)$ is a Lorentzian manifold with d'Alembertian
$\Box$, then the symbol of $\Box$ is the metric tensor in the contravariant
form. In this case the bicharacteristic segments are the null geodesic
segments, and the inhomogeneous wave equation $\Box u = f$ has global
solutions in the distribution sense if (1) $(M,g)$ is null pseudoconvex
and (2) each end of each inextendible null geodesic fails to be
imprisoned. This second requirement is called {\em disprisonment\/} of
null geodesics. In the language of PDE's, it is the requirement that the
operator be of real principal type. Beem and Parker \cite{KGSGG,GBSS,WSS} have
obtained several powerful results on the stability of solvability of
pseudodifferential
equations {\em via\/} stability theorems for pseudoconvexity and
disprisonment. This approach also yielded some new methods in the study of
sectional curvature \cite{VPSC}.

Interestingly, both pseudoconvexity and disprisonment of null geodesics
fail to be separately $C^1$-stable in the Whitney topology, but the
requirement that they hold jointly {\em is\/}  $C^1$-stable \cite{WSS}.

As one would expect, disprisonment and pseudoconvexity have important
implications for the geodesic structure of a spacetime. Williams \cite {W}
found examples of geodesically complete spacetimes with arbitrarily close
incomplete metrics in the Whitney $C^{r}$-topology. Hence geodesic
completeness fails to be Whitney $C^{r}$-stable for all $r\geq 1$. On
the other hand, Beem and Ehrlich \cite{BE} have established the
$C^{1}$-fine stability of causal geodesic completeness for Lorentzian
manifolds which are both causally pseudoconvex and causally disprisoning.

Beem and Parker \cite{PGR} established an extension of
Seifert's result \cite{S} to manifolds with a linear connection which is
pseudoconvex and disprisoning for {\em all\/} types of geodesics. In
particular, one now has pseudoriemannian versions of the Hopf-Rinow and
Hadamard-Cartan Theorems.

\subsection{\subheads Sprays}
One of the most important generalizations of
ordinary differential equations to manifolds is the well-known class of
second-order differential equations or sprays. For example, they
occur as the Hamiltonian vector fields of regular Lagrangians in
variational problems \cite{AB,K}. A (general) {\em spray\/} on a manifold
$M$ is
defined as a projectable section of the second-order tangent bundle
$TTM \rightarrow TM$. This is precisely the condition needed to define a
second-order differential equation \cite{BC,AB}. Recall that an integral curve
of a
vector field on $TM$ is the canonical lift of its projection if and only
if the vector field is projectable. For
any curve $c$ in $M$ with tangent vector field $\dot c$, this $\dot c$ is
the canonical lift of $c$ to $TM$ and $\ddot c$ is the canonical lift of
$\dot c$ to $TTM$. Then each projectable vector field $S$ on $TM$
determines a second-order differential equation on $M$ by $\ddot c  =
S\circ \dot c$ and any such curve
with $\dot c(s_{0}) = v_{0}\in T_{c(s_{0})}M$
is a solution with initial condition $v_0$. Solutions are preserved
under translations of parameter, they exist for all
initial conditions by the Cauchy theorem, and, as our manifolds are
assumed to be Hausdorff, each solution will be unique provided we
take it to have maximal domain; {\em i.e.,} to be inextendible
\cite{BC,DR,DRD}.

Let $J$ be the canonical involution
on $TTM$ and $C$ the Euler  (or Liouville) vector field on $TM$. We recall
that in local coordinates, $J(x,y,X,Y) = (x,X,y,Y)$ and
$C:(x,y)\mapsto (x,y,0,y)$. Then a section $S$ of $TTM$ over $TM$ is
a {\em spray} when $JS  =  S$; that is, when it can be
expressed locally as $S:(x,y)\mapsto (x,y,y,Y(x,y))$.
 We say that a spray $S$ is (positively) {\em homogeneous}
of degree $m$ when $[C,S] = (m-1)S$ for $m \ge 0$.
In this case the functions $Y(x,y)$ are homogeneous of degree
$m$ in the fiber component: $Y(x,ay)  =  a_{*}a^{m-1}Y(x,y)$
\cite {KV}. Here $a_{*}$  denotes the induced tangent map of scalar
multiplication by $a$. We denote the set of sprays on $M$ by $\sp(M)$ and
those homogeneous of degree $m$ by $\sp_m(M)$.
It has been usual to consider only
(positive) integral degrees of homogeneity, but we make no such
restriction. Previously \cite{DR,KV}, our
(general) sprays have been called {\em semi\-sprays\/} and the name
{\em sprays\/} reserved for those homogeneous
of degree two. We do not make this restriction either. We do, however,
consider only sprays defined on the entire tangent bundle $TM$; others
\cite{KV} have used
the reduced tangent bundle with the 0-section removed.

For some purposes, it is more convenient to use a different
characterization of sprays \cite{DR,KV}. The vertical endomorphism $V$,
in local coordinates given by $V(x,y,X,Y)  =  (x,y,0,X)$, may be
regarded as a vector-valued 1-form on $TM$. We observe that $V:TTM
\rightarrow {\scr V}TM$, the vertical bundle over $TM$, and is a
nilpotent map: $V^{2} = 0$. Then a spray also can be characterized by
$VS = C$. This version has been used, for example, in stability
theory \cite{DRD}.

Recall that in general a connection only provides a horizontal subbundle
of $TTM$ complementary to the vertical subbundle.
The Nijenhuis bracket ({\em e.g.,} \cite{DRD}) determines for each spray
(or connection) a Lie
subalgebra of the Lie algebra of vector fields on $TM$. This subalgebra
consists precisely of those morphisms of $TTM$ over $TM$ which preserve
the horizontal and vertical subbundles \cite{DR}.

Several important results concerning sprays \cite{APS,BC,D,KV} rely on
the facts that each spray $S$
determines a unique torsion-free connection $\Gamma$, and conversely, every
spray $S$ arises from a connection $\Gamma$ the torsion of which can be
assigned arbitrarily. The solution curves of the differential equation
$\ddot c  =  S_{\Gamma}\circ \dot c$ for a connection-induced spray are
precisely the geodesics of that connection. The familiar geodesic spray
corresponding to a linear connection is a quadratic spray: $[C,S] = S$.
In this case its solution curves are not only preserved under
translations, but also under affine transformations of the parameter
$s \mapsto as+b $ for constants $a,b$ with $a\neq 0$.

Here is perhaps the simplest example of a spray arising from a regular
variational problem.
If $(M,g)$ is a pseudoriemannian manifold, then the energy function
$\epsilon _{g}$ is defined by
$$\epsilon _{g}: TM \rightarrow \R: v
\mapsto\textstyle\frac{1}{2} g(v,v)\, .$$
The canonical spray $S_g$ on $M$ is defined as the vector field on $TM$
corresponding
to the 1-form $-d\epsilon_{g}$ on $T^{*}M$ with respect to the
canonical symplectic structure on $T^*M$. As a derivation
on real functions defined on $TM$, $S_{g}$ annihilates $\epsilon_g$.
Thus the lifts $\dot c$ of solution
curves $c$ are
integral curves in $TM$ of $S_{g}$ along which $\epsilon _{g}$ is constant.
Now the Levi-Civit\`a connection
$\Gamma_{g}$ determines a unique spray which also annihilates $\epsilon_g$.
It follows that $S_g$ is the geodesic spray, so the solution curves $c$
for $S_{g}$ are the geodesics of $\Gamma_{g}$.

\bigskip\medskip\noindent
In this paper, we shall study the combination: general homogeneous
sprays which are pseudoconvex and disprisoning. Among other things, this
may be regarded as a continuation of the program to geometrize the study
of PDE's begun by Beem and Parker \cite{KGSGG,GBSS,WSS,PGR}.
Section 2 contains our
definitions, notations and conventions. Section 3 is devoted to the
generalization of the main results of \cite{PGC} to a large class
of sprays. Finally, Section 4 gives a generalization of the main
stability result of \cite{WSS} to homogeneous sprays.

Throughout, all manifolds are smooth (meaning $C^{\infty}$), connected,
paracompact, Hausdorff, and usually noncompact (see Section 2).

The authors would like to thank Boeing Wichita and the National Science
Foundation
for travel and support grants, and The Wichita State University, CICY,
CIMAT, and Universidad Aut\'onoma de San Luis Potos\'{\i} for hospitality
during the progress of this work.

\section{\heads Preliminaries}

We begin with the principal definitions. Let $S$ be a spray on $M$.
\begin{definition}
We say that a curve $c:(a,b)\rightarrow M$ is a\/ {\em geodesic} of $S$ or
an {\em $S$-geodesic} if and only if the natural lifting $\dot c$ of $c$ to
$TM$ is an integral curve of $S$.
\end{definition}
This means that if $\ddot c$ is the natural lifting of $\dot c$ to $TTM$,
then $\ddot c = S(\dot c)$.
\begin{definition}
We shall say that $S$ is\/ {\em pseudoconvex} if and
only if for each compact $K \subseteq  M$ there exists a compact $K^\prime
\subseteq  M$
such that each $S$-geodesic segment with both endpoints in $K$ lies
entirely within $K^{\prime}$.
\end{definition}
If we wish to work directly with the integral curves of $S$, we merely
replace ``in''  and ``within'' by ``over''.
\begin{definition}
We shall say that $S$ is\/ {\em disprisoning} if and
only if no inextendible $S$-geodesic is contained in (or lies over) a
compact set of $M$.
\end{definition}
In relativity theory \cite{HE}, such inextendible geodesics are said to be
imprisoned in compact sets; hence our name for the negation of this property.

Following this definition, we make a convention: all $S$-geodesics are to
be regarded as always extended to the maximal parameter intervals
({\em i.e.,} to be
inextendible) unless specifically noted otherwise.
When the spray $S$ is clear from context, we refer simply to geodesics.
Also, we shall consider only noncompact manifolds because no
spray can be disprisoning on a compact manifold. However,
Corollary~\ref{qcov}
may be used to obtain results about compact manifolds for which the
universal covering is noncompact.
\begin{example}
{\ns When $S$ is a quadratic spray, we recover the
notions
previously defined by Beem and Parker \cite{PGC} for linear
connections.}
\end{example}

The natural lift $\dot c$ of a curve $c$ was denoted by $c'$ in
\cite{PGC}. With this change in notation, the proof of Lemma 3 there
applies word-for-word to sprays. Thus we have
\begin{lemma}
Let $M$ be a manifold with a pseudoconvex and disprisoning
spray $S$. If $p \neq q$ , $p_n \rightarrow p$ and
$q_n \rightarrow q$, and if for each $n$ there is a geodesic segment
from $p_n$ to $q_n$, then there is a geodesic segment from $p$ to
$q$.\eop
\label{l3}
\end{lemma}
As in \cite{KGSGG}, we obtain
\begin{proposition}
Let $M$ be a manifold with a pseudoconvex and disprisoning spray $S$.
If $K \subseteq  M$ is compact, then the geodesic convex hull
${[\! [}  K{]\! ]} $ is compact.
\end{proposition}
\begin{proof}
Pseudoconvexity implies that ${[\! [}  K{]\! ]} $ is contained in a compact
set,
and the lemma shows ${[\! [}  K{]\! ]}$ is closed, hence compact.
\end{proof}

\section{\heads Geodesic Systems}

In \cite{PGC}, results were obtained concerning geodesic connectedness of
manifolds with a linear connection. Since quadratic sprays are equivalent
to linear connections, we immediately obtain corresponding results
for quadratic sprays. For example \cite[Proposition 5]{PGC}
\begin{proposition}
Let S be a pseudoconvex and disprisoning quadratic spray on $M$. If $S$
has no conjugate points, then $M$ is geodesically connected. In other
words, every pair of points in $M$ may be joined by at least one
$S$--geodesic segment.\eop
\label{qncpgc}
\end{proposition}
However, we are primarily interested in more general sprays.

Recall that for each $p\in M$ and each $v\in T_pM$, there is a unique
geodesic $c_v$ such that $c_v(0) = p$ and $\dot{c}_v(0) = v$.
\begin{definition}
The {\em exponential map} of $S$ at $p$ is given by $\exp_p(v) = c_v(1)$
for all $v\in T_pM$ such that $c_v(1)$ exists.
\end{definition}
Thus, as in the usual cases, the domain of $\exp$ is an open tubular
neighborhood of the 0-section in $TM$. The next lemma follows from
Lemma~\ref{l3} as in \cite{PGC}.
\begin{lemma}
Let $M$ be a manifold with a pseudoconvex and disprisoning
spray $S$. Assume $p \neq q$ and $q_n \rightarrow q$. If\/ $(v_n)$ is a
sequence
in $T_p M$ such that $\exp_p(v_n) = q_n$, then there is a vector $v\in T_p M$
and a subsequence $(v_k)$ of $(v_n)$ such that $v_k\rightarrow v$ and
$\exp_p(v) = q$.\eop
\label{l4}
\end{lemma}
\begin{theorem}
If $S$ is a homogeneous spray, then its exponential map $\exp$ is a
local diffeomorphism.
\end{theorem}
\begin{proof}
We proceed as in the usual proof ({\em e.g.,}
\cite[p.\,116]{BJ}) except that now $\exp_p(tv) = c_{tv}(1) = c_v(t^m)$
for $0\le t\le 1$ when $S$ is homogeneous of degree $m$. But we still
obtain $\exp_{p*} = 1$ as usual.
\end{proof}
All examples of general sprays which we have examined have this property.
We conjecture that it is true of all sprays, but we cannot
prove it yet. Thus we make
\begin{definition}
A spray is\/ {\em LD} if and only if its exponential map is a local
diffeomorphism.
\label{ld}
\end{definition}
What we shall use is that the goeodesics of such sprays give normal starlike
neighborhoods of each point in $M$. This fact together with Lemma~\ref{l3}
yields the next result, as in \cite[Prop.\,5]{PGC}.
\begin{proposition}
Let $M$ be a manifold with a pseudoconvex and disprisoning LD spray $S$. If\/
$S$ has no conjugate points, then $M$ is geodesically connected.\eop
\label{p5}
\end{proposition}

Let $M$ be a manifold with a spray $S$ and let
$\widetilde {M}$ be a covering
manifold. If $\phi:\widetilde{M}\rightarrow M$ is the covering map, then it
is a local diffeomorphism. Thus $\tilde{S} = (\phi_{\ast})^{\ast}S$ is
the unique spray on $\widetilde{M}$ which covers $S$, geodesics of
$\tilde{S}$ project to geodesics of $S$ and geodesics of $S$ lift to
geodesics of $\tilde{S}$. Also, $S$ has no conjugate points if and only
if $\tilde{S}$ has none. The fundamental group is simpler, and $\tilde{S}$
may be both pseudoconvex and disprisoning even if $S$ is neither.
Proposition~\ref{p5} and simple projection arguments yield
({\em cf.} \cite[Corollary 6]{PGC})
\begin{corollary}
Let $M$ be a manifold with a pseudoconvex and disprisoning
LD spray $S$ and let $\widetilde{M}$ be a covering
manifold with covering homogeneous spray $\tilde{S}$. If
$\tilde S$ has no conjugate
points, then both $\widetilde M$ and $M$ are geodesically connected.\eop
\label{qcov}
\end{corollary}
The next result is the analogue of Theorem 9 of \cite{PGC}.
\begin{theorem}
Let $S$ be a pseudoconvex and disprisoning LD spray on $M$. If
$S$ has no conjugate points, then for each $p \in M$ the exponential map
of $S$ at $p$ is a diffeomorphism.\eop
\end{theorem}
We remark that none of these results require (geodesic) completeness of
the spray $S$.

\section{\heads Stability}

In this section we consider the joint stability of pseudoconvexity and
disprisonment for homogeneous sprays in the fine topology. Because each
linear
connection determines a homogeneous spray, Examples 2.1 and 2.2 of
\cite{WSS} show that
neither condition is separately stable. (Although \cite{WSS} is written
in terms of principal symbols of pseudodifferential operators, the cited
examples are actually metric tensors).  We shall obtain
$C^{0}$-fine stability, rather than $C^{1}$-fine stability as in \cite
{WSS}, due to our effective shift from potentials to fields as the basic
objects. The proof requires only minor modifications of that in \cite
{WSS}, so we shall concentrate on the changes here and refer to \cite
{WSS} for an outline and additional details.

Rather than considering $r$-jets of functions, we now take $r$-jets of
sections in defining the Whitney or $C^r$-fine topology
as in Section 2 of \cite{WSS}. Also, just as was done there, we must
modify these topologies to take homogeneity into account. Let $h$ be an
auxiliary complete Riemannian metric on $M$. A homogeneous
spray is determined by its degree of homogeneity $m$ and its restriction to
the $h$-unit sphere bundle $UM$ in $TM$. (Note that our $UM$ replaces
$S^*M$ in \cite{WSS}, changing from the cotangent to the tangent bundle.)
Thus we actually look at the
$C^r$-fine topology on the sections of $TTM|UM$ over $UM$. However, as in
\cite{WSS}, we shall say that a set in $\sp_m(M)$ is open if and only if the
corresponding set in the sections over $UM$ is open.

If $\g_1$ and $\g_2$ are two integral curves of a spray $S$ with $\g_1(0)
= (x,v)$ and $\g_2(0) = (x,\l v)$ for some positive constant $\l$, then
the inextendible geodesics $\p\circ\g_1$ and $\p\circ\g_2$ differ only by
a reparametrization. Thus, as in \cite{WSS}, it will suffice to consider
only one integral curve for each direction at each point of $M$.

Observe the the equations of geodesics involve no derivatives of $S$. Thus
if $\g :[0,a]\rightarrow TM$ is a fixed integral curve of $S$ in $TM$ with
$\g(0) = v_0 \in UM$ and if $\g' :[0,a]\rightarrow TM$ is an integral
curve of $S'$ in $TM$ with $\g'(0) = v$, then $d_h \left( \p\circ\g(t),
\p\circ\g'(t)\right) < 1$ for $0\le t\le a$ provided that $v$ is
sufficiently close to $v_0$ and $S'$ is sufficiently close to $S$ in the
$C^0$-fine topology. This and the compactness of $UK_1$ when $K_1$ is
compact yield the following result.
\begin{lemma}
Assume $K_{1}$ is a compact set contained in the interior of the
compact set $K_{2}$ and let $S$ be a disprisoning homogeneous spray.
There exist
tangent vectors $v_{1},...,v_{m} \in TK_{1}$ and positive constants
$\delta_{1},...,\delta_{m},\alpha_{1},...,\alpha_{m},\epsilon$ such that
if $S'$ is in a $C^{0}$-fine $\epsilon$-neighborhood of $S$ over
$V$, then the following hold:
\begin{enumerate}
\item if $c$ is an inextendible $S$-geodesic with $c(0)$
in a
$\delta_{i}$-neighborhood of $v_{i}$, then $c[0,a_{i}] \subset V$
and $c (a_{i}) \in V-K_{2}$;
\item If $c^{\prime}$ is an inextendible $S^{\prime}$-geodesic with
$\dot{c}^{\prime}(0)$ in a $\delta _{i}$-neighborhood if $v_{i}$, then
$c^{\prime}[0,a_{i}] \subset V $ and $c^{\prime}(a_{i}) \in V-K_{2}$;
\item Two inextendible geodesics, $c$ of $S$ and $c^{\prime}$ of
$S^{\prime}$ with $\dot{c}(0)$ and $\dot{c}'(0)$ in a $\delta
_{i}$-neighborhood of $v_{i}$, remain uniformly close together for
$0\leq t \leq a_{i}$;
\item The union of all the $\delta_{i}$-neighborhoods of the $v_{i}$ is
large enough to cover the part of $TK_{1}$ in which we are interested.\eop
\end{enumerate}
\end{lemma}

Continuing to follow \cite{WSS}, we construct the increasing sequence
of compact sets $\{ A_n\}$ which exhausts $M$ and the monotonically
nonincreasing sequence of positive constants $\{ \epsilon_n\}$. The only
changes from \cite[p.\,17f\,]{WSS} are to use
integral curves of $S$ in $TM$ instead of bicharacteristic strips in $T^*M$.
No additional changes are required for the proof of the next result either.
\begin{lemma}
Let $S$ be a pseudoconvex and disprisoning homogeneous spray and let
$S^{\prime}$
be $\delta$-near to $S$ on $M$. If $c^{\prime}:(a,b)\rightarrow M$ is an
inextendible $S'$-geodesic, then there do not exist values $a <
t_{1} <
t_{2} < t_{3} < b$ with $c^{\prime} (t_{1}) \in A_{n}$, $c ^{\prime}
(t_{3}) \in A_{n}$ , and $c^{\prime} (t_{2}) \in A_{n+4} - A_{n+3}$.\eop
\label{3.2}
\end{lemma}

Now we establish the stability of pseudoconvex and disprisoning
homogeneous sprays by showing that the set of all sprays in $\sp_m(M)$
which is pseudoconvex and disprisoning is an open set in the $C^0$-fine
topology. The only changes needed from the proof of Theorem 3.3 in
\cite[p.\,19]{WSS} are replacing principal symbols by sprays,
bicharacteristic strips by integral curves, $S^*A_n$ by $UA_n$, and
references to Lemma 3.2 there by references to Lemma~\ref{3.2} here.
\begin{theorem}
If $S\in\sp_m(M)$ is a pseudoconvex and disprisoning homogeneous spray,
then there is some $C^{0}$-fine neighborhood $W(S)$ in $\sp_m(M)$ such
that each $S' \in W(S)$ is both pseudoconvex and disprisoning.\eop
\label{spd}
\end{theorem}

Since linear connections may be identified with quadratic sprays, we
immediately obtain
\begin{corollary}
Pseudoconvexity and disprisonment are jointly $C^0$-fine stable properties
in the space of linear connections.\eop
\end{corollary}
In particular,
\begin{corollary}
If $M$ is a pseudoconvex and disprisoning pseudoriemannian manifold,
then any linear connection on $M$ which is sufficiently close to
the Levi-Civit\`a connection is also pseudoconvex and disprisoning.\eop
\end{corollary}

If we denote by $\sp_{H}(M)$ the set of {\em all\/} homogeneous sprays
of any degree $m$ on $M$, then we may topologize it by taking the weak topology
generated by those on the subsets $\sp_m(M)$. Now Theorem~\ref{spd} can be
given a somewhat more general formulation.
\begin{theorem}
If $S$ is a pseudoconvex and disprisoning homogeneous spray, then any
sufficiently close homogeneous spray (of\/ {\em any} degree of homogeneity)
is also pseudoconvex and disprisoning.\eop
\end{theorem}
\begin{corollary}
If $(M,g)$ is a pseudoconvex and disprisoning pseudoriemannian manifold,
then any homogeneous spray $S$ which is sufficiently close to the geodesic
spray $S_g$ is also pseudoconvex and disprisoning.\eop
\end{corollary}

\end{document}